\documentclass[11pt, a4paper]{article}
\usepackage[utf8]{inputenc}
\usepackage{amsmath,amssymb,amsfonts,dsfont}
\usepackage{amsthm}
\usepackage{mathtools}
\usepackage{verbatim}
\usepackage{a4wide}
\usepackage{epsfig,epic,eepic,graphicx}
\usepackage[mathcal]{euscript}
\usepackage{import}
\usepackage{xifthen}
\usepackage{comment}
\usepackage[usenames,dvipsnames]{xcolor}
\usepackage[unicode, bookmarks, colorlinks, breaklinks, pagebackref,
]{hyperref}  
\hypersetup{linkcolor=blue,citecolor=blue,filecolor=black,urlcolor=Magenta}

\usepackage[capitalize,nameinlink,noabbrev]{cleveref}

\usepackage{cite}

\usepackage{enumitem}

\newtheorem{theorem}{Theorem}[section]
\newtheorem{lemma}{Lemma}[section]

\newtheorem{definition}{Definition}

\numberwithin{equation}{section}
\makeatletter

\@addtoreset{equation}{section}
\makeatother
\newcommand{\email}[1]{\href{mailto:#1}{\texttt{#1}}}

\def\XXint#1#2#3{{\setbox0=\hbox{$#1{#2#3}{\int}$ }
\vcenter{\hbox{$#2#3$ }}\kern-.6\wd0}}

\DeclarePairedDelimiterX\Set[1]\{\}{%
  #1%
}
\allowdisplaybreaks

\begin{document}
\title{Bandgaps at high frequency from highly oscillatory metamaterials}
%\title{Band gaps from highly oscillatory, high frequency metamaterials}

 \author{Silvia Jim\'{e}nez Bola\~{n}os\footnote{Department of Mathematics,
    Colgate University, 13 Oak Drive, Hamilton, New York 13346, USA (\email{sjimenez@colgate.edu}). } 
\and
 Robert P. Lipton\footnote{Department of Mathematics, Louisiana State University, 265 Lockett Hall
Baton Rouge, Louisiana 70803, USA
    (\email{lipton@lsu.edu}).}
  \and
  Robert Viator\footnote{Department of Mathematics, Denison University, 100 West College Street, Granville, Ohio 43023, USA (\email{viatorr@denison.edu}). }}

\date{}

\maketitle

\begin{abstract}
%  We consider the Bloch spectrum for TE-modes in a 2D high-constrast dielectric metamaterial, in the case where the propagating wave has wavelength commensurate with the length scale $d$ of the metamaterial microstructure and the permittivity in the inclusions scales as $d^{-2}$.  We show that the Bloch eigenvalues can be expanded in a power series in $d^{2}$, provided the Bloch quasi-momentum satisfies an upper bound given in terms of the inclusion geometry.  We also prove that band gaps appear in the spectrum for sufficiently small $d<\overline{d}$, where $\overline{d}>0$ is prescribed in terms of the Dirichlet and Neumann-Poincar\'e spectrum of the inclusion domain.

We investigate the Bloch spectrum for transverse electric (TE) modes in a 2D high-contrast dielectric metamaterial. Our analysis focuses on the high-frequency regime where the propagating wavelength is commensurate with the microstructural period $d$, and the permittivity within the inclusions scales as $d^{-2}$. Using operator factorization and Sobolev space decomposition, we demonstrate that Bloch eigenvalues $\omega^2$ can be analytically continued and represented as a power series in $d^2$. We establish a rigorous upper bound on the microstructural length scale $d$, defined by the inclusion geometry, under which this expansion remains valid. Furthermore, we provide a definitive proof for the opening of band gaps in the frequency spectrum for all periods below a critical length scale, i.e.,
 for all $d < \overline{d}$. The threshold is explicitly characterized as a geometric function determined by the interaction between Dirichlet eigenvalues and the spectrum of the associated Neumann-Poincaré operator.
\end{abstract}
\section*{Acknowledgement}
This material is based upon work supported in part by the National Science Foundation under Grant Numbers DMREF-1921707 and DMS-2110036.

\section{Introduction}
\label{sec:introduction}

Metamaterials have seen an increase in use and implementation in  science, mathematics, and engineering disciplines, ranging from acoustics \cite{Smyshlyaev2009, Smith2011} to electromagnetism \cite{Pendry2000,Kadic2019,Bouchitte2010} to heat transport and mechanics \cite{Ying2021,Guenneau2011,Davies2025,Kamotsli2019,Guenneau2012}.  A defining characteristic of these materials is their exotic ``averaged" or effective macroscopic response, usually represented by using an effective operator derived using multiscale homogenization methods, which arises from highly oscillatory heterogeneities in the metamaterial at very fine length scales with respect to the phenomenon at hand.  %Notable examples include negative refractive index materials in acoustics and electromagnetism \textcolor{red}{(again, Pendry and others \cite{Pendry2000}, include the Craster \& Guenneau book)}, and thermal metamaterials designed to manipulate heat flux through highly contrasted conductivity tensors \textcolor{red}{(citations, yo!)}. 

In this paper, we analyze a two-dimensional, high-contrast electromagnetic metamaterial.  However, diverging from the classical sub-wavelength assumptions, %found (again, cite Pendry, others \cite{Pendry2000}), 
we investigate the scattering regime where the wavelength is \emph{commensurate} with the length scale of the microstructure (the wavelength $\lambda$ is commensurate to the microstructural period $d$, i.e. $\lambda \approx d$). This regime has previously been explored with respect to its effect on higher length scales through \emph{High-Frequency Homogenization} (HFH), a technique initiated by Craster et al. See \cite{craster2010,craster2016,Assier2020,touboul2023highfrequency} as well as very long time asymptotic behavior of the wave equation with highly oscillatory initial data \cite{APR2008}. In \cite{craster2010}, the authors perturbed from standing waves supported by the crystal at wavelengths commensurate with the length scale of the crystal; the asymptotics developed there were used to accurately reproduce the Bloch modes near the edges of the Brillouin zone, which allows for the homogenization of periodic media in the vicinity of ``cell resonances".  This method was extended to the scalar and vector wave equations and Schrodinger-type equations in \cite{craster2016}, to acoustic spring-mass models with imperfect interfaces in \cite{Assier2020}, and to frequency-dependent dispersive media in \cite{touboul2023highfrequency}.

The problem addressed here is different; we seek to understand the spectral properties of metamaterials across higher frequencies as an operator perturbation problem,  \textcolor{blue}{without multiscale analysis.} In particular, in this paper we assert the %possibility of opening up 
existence of band gaps over higher frequency intervals in the spectrum. Here, we present the model for the two-dimensional metamaterial, and we complexify the length scale as well as the material properties to arrive at a power series representation of the solution. Returning to real periods and rescaling, we arrive at a spectral problem that is strikingly similar to those seen in \cite{RobertRobert1, RobertRobert2}, albeit with coefficients that now explicitly depend on the quasimomentum of the Bloch waves.  Following the methodology in \cite{RobertRobert1, RobertRobert2}, we develop an analytic power series expansion for Bloch modes relative to the high-contrast parameter. A key departure from previous studies is that this radius of convergence allows us to obtain an upper bound on wavelengths for which the power series expansion for the Bloch modes will converge. The bound is explicit and stated in terms of the geometry of the metamaterial.  For all wavelengths within this bound, we have a uniform bound on the error incurred approximating by finite sums of terms in the power series expansion of the Bloch modes.  We then leverage the power series convergence and corresponding error analysis, together with a spectral problem which serves as an ``upper bound" to the Bloch spectrum, to describe gaps in the spectrum for the metamaterial.

Central to prescribing both the wavelength bound and the spectral gaps are the \emph{electrostatic resonances of the periodic structure}, i.e. functions $v$ which satisfy either a quasi-periodic or homogeneous Neumann boundary condition on the boundary of each period cell, and real eigenvalues $\mu$ for which
\begin{equation}\label{electrostatic}
-\nabla \cdot (\chi_{D} \nabla v) = -\mu \Delta v \text{,}    
\end{equation}
where $D$ is a disconnected component of the metamaterial in a given period cell. On expressing \eqref{electrostatic} in terms of double layer potentials shows resonances are connected to the spectra of Neumann-Poincar\'e operators associated with double layer potentials as given in \cite{kanglayerpotential13} and \cite{khavinsonputinarshapiro07}.  The are the well-known electrostatic resonances used in the composites literature and used to bound effective properties; see c.f. \cite{bergmancomp78},\cite{bergmancube79},\cite{brunocomp91}, \cite{goldenpapa83}, \cite{mcphedranmilton81}, \cite{miltontheorycomposites02}.
%The exterior is modeled as an ENZ material surrounding islands of unit dielectric constant. The $z$ will be a ratio of dialectric constants.  

%Need to discuss ENZ materials \cite{KohnVenkatraman2024,kohn2024transverse} \textcolor{blue}{are we still working in this context?  It looks like (from eq. (2.1)) our permittivity here is quite large.}

%We need to discuss high-frequency homogenization \cite{craster2010,craster2016,Assier2020,touboul2023highfrequency} \textcolor{blue}{This is done, RV}
%\subsection{Outline of the paper}

Our main result, found in section \ref{s:Bandgaps}, is the presence of band gaps in the frequency spectrum of the metamaterials presented in the paper, along with a bound on the length scaling $d$ necessary to obtain these spectral gaps.  The gaps are present for a wide class of inclusion geometries, and their size and location are estimated by the Neumann-Poincar\'e and Dirichlet spectra of the inclusions.  We present the main theorem specialized to disk-shaped inclusions in the next section to illustrate the ideas.   

The paper is organized as follows. In section \ref{sec:formulation}, we introduce the mathematical formulation of the ``high-frequency" analysis of the spectrum of transverse-electric modes propagating through the metamaterial, rescaling the problem to obtain a sesquilinear form with $\mathbb{Z}^2$-periodic coefficients on $L^2(\mathbb{R}^2)$ and we provide an application of the main results to the particular case of buffered disks inclusions.  In section \ref{s:PowerSeriesRep}, we invoke the analysis of previous work \cite{RobertRobert1}, \cite{RobertRobert2} to establish power series expansions for the Bloch eigenmodes associated with the form from section \ref{sec:formulation}, complete with a radius of convergence for each quasi-momentum $\alpha$:  the limit spectra for these power series are obtained in section \ref{s:limitspectra}.  In section \ref{s:wavelengthbound}, we obtain a lower bound on the wavelength (and thus an upper bound on the frequency) of incident waves for which the analysis is valid - this is natural in light of previous work on high-frequency homogenization (\cite{Assier2020}, \cite{Guenneau2011},\cite{craster2010}, \cite{craster2016}).  Finally, in section \ref{s:Bandgaps}, we prove the existence of bandgaps for a broad class of inclusion geometries, presenting finer details in the case of disk-shaped inclusions in section \ref{sec:formulation}.

\section{Problem Formulation and Complexification of the Period.}
\label{sec:formulation}

Band-gap behavior in the metamaterial crystal arises from multiple scattering effects. To capture these phenomena, we consider a regime where the lattice period $d$ is commensurate with the wavelength of the incident light.

We begin by describing the two-dimensional periodic lattice transverse to parallel rods on the unit cell $Y=[0,1)^2$, and subsequently rescale. Let $D\subset \subset Y$ be open and bounded with $C^{1,\gamma}$ boundary for $\gamma>0$, denoting the rod cross-section scaled to the unit period, and let $D_d\subset Y_d$ represent the corresponding $d$-periodic cross-section. We define $\Omega$ as the union of all cross sections of the rods in the unit-periodic lattice
\begin{equation*}    \Omega=\bigcup_{\boldsymbol{k}\in\mathbb{Z}^2}D + \boldsymbol{k},
\end{equation*}
and $\Omega_d$ as the union of all cross-sections in the $d$-periodic lattice 
\begin{equation*}    \Omega_d=d\Omega=\bigcup_{\boldsymbol{k}\in\mathbb{Z}^2}d(D+\boldsymbol{k}).
\end{equation*}

We characterize the band gaps and pass bands associated with this inclusion geometry under transverse electric (TE) polarized electromagnetic radiation where the magnetic field is oriented parallel to the rods and the electric field lies in the plane transverse to the rods.

The connected region separating the rods, namely $\mathbb{R}^2\setminus\Omega_d$, is referred to as the matrix.  %Denote by $D_d=d\,D$ and $Y_d=d\,Y$.  
The direction of propagation of the TE wave in the transverse  plane is described by the unit vector $\boldsymbol{\kappa}\in{\mathbb{R}}^2$. We denote the wave number by $k=2\pi/\lambda$, where $\lambda$ is the wavelength.

The relative dielectric permittivity, denoted by $\epsilon_d$, depends on the lattice period $d$ and is given by
\begin{equation} 
\label{dieconsa}
\epsilon_d^{-1} = 
     \begin{cases}
       d^2\eta &\quad\text{in}\quad D_d\\
       1 &\quad\text{in}\quad Y_d\setminus D_d,
     \end{cases}
\end{equation}
with $\eta\in\mathbb{C}$.  For TE modes, the electric and magnetic fields reduce to $E_d=(E_1,E_2,0)$ and $H_d=(0,0,v_d)$, where the scalar function $v_d$ satisfies the Helmholtz equation 
\begin{equation}
	\label{Helm-d}
	-\nabla\cdot\left(\epsilon_d^{-1}(x/d)\nabla v_d(x)\right)=\omega^2 v_d(x), \hspace{2mm}x\in\mathbb{R}^2.
\end{equation}
This spectral problem describes wave propagation inside the crystal at frequency $\omega$. The self-adjoint divergence-form operator associated with \eqref{Helm-d} is 

\begin{equation}
\tilde{L}_d = -\nabla \cdot ((1 - \chi_{\Omega_d}) + d^2\eta\, \chi_{\Omega_d})\nabla \label{tildaLd}
\end{equation}
where $\chi_A$ represents the characteristic function of the set $A$; i.e, $\chi_A(x)=1$ if $x\in A$ and $\chi_A(x)=0$ otherwise.
%Add the $\kappa\tau$ here. See Chapter 11 in Milton's book page 223.
By Bloch's Theorem, the magnetic field is quasi-periodic and can be written as $$v_d(x)=
\varphi_d(x)e^{ik\boldsymbol{\kappa}\cdot x},$$ where $
    \varphi_d$ is $d$-periodic.  Introducing the change of variables $y=\frac{x}{d}$, we have $\partial_{x_i}=\frac{1}{d}\partial_{y_i}$, for $i=1,2$, and equation \eqref{Helm-d} becomes
\begin{equation}
	\label{Helm-d2}
	-d^{-2}\nabla_y\cdot\left(\epsilon_d^{-1}(y)\nabla_y v_d(d\,y)\right)=\omega^2 v_d(d\,y).
\end{equation}
We define the rescaled fields $u_d(y):=v_d(d\,y)$ and $\psi_d(y):=\varphi_d(d\,y)$ so that
\begin{equation*}
u_d(y)=\psi_d(y)e^{i\,dk\boldsymbol{\kappa}\cdot y}.
\end{equation*}
The rescaled field $u_d$ satisfies the equation on the unit cell $Y$,
    \begin{equation}
	\label{Helm-d3}
	-\nabla_y\cdot\left(a_d(y)\nabla_y u_d(y)\right)=\omega^2 u_d(y),
\end{equation}
with coefficients
\begin{equation*}  
a_d(y) = 
     \begin{cases}
       \eta &\quad\text{in}\quad D\\
       \frac{1}{d^{2}} &\quad\text{in}\quad Y\setminus D.
     \end{cases}
\end{equation*}

This rescaled formulation allows us to take advantage of the analyticity of the operator $-\nabla_y\cdot a_d\nabla_y $ acting on quasi-periodic functions, as described in \cite{RobertRobert1,RobertRobert2}.
To proceed further, we complexify the parameter $d$ by introducing $z=d^{2}\in\mathbb{C}$.  We then extend the quasi-momentum vector by setting $\alpha:=dk\boldsymbol{\kappa}=\sqrt{z}k\boldsymbol{\kappa}$ so that, for real positive $z$, we have that $\alpha$ lies in the first Brillouin zone of the reciprocal lattice given by $\in Y^* =(-\pi,\pi]^2$. 

Setting $h=1/z$, we have $\alpha=\frac{k\boldsymbol{\kappa}}{\sqrt{h}}$. The self-adjoint divergence-form operator associated with \eqref{Helm-d2} is 
\begin{equation}\label{OpeLh}
    L_h=-\nabla\cdot(\eta\chi_{\Omega}+h(1-\chi_{\Omega}))\nabla 
\end{equation}
defined by the sesquilinear form over $L^2(\mathbb{R}^2)\times L^2(\mathbb{R}^2)$ given by
\begin{equation}
 B_{\eta,h}(u,w)=   \int_{\mathbb{R}^2}(\eta\chi_{\Omega}+h(1-\chi_{\Omega})\nabla u(x)\cdot\nabla\overline{w}(x)\,dx
 \label{Betah}
\end{equation}
with domain $W^{1,2}(\mathbb{R}^2)$. 

Floquet theory shows that the spectrum $\sigma(L_h)$ has the band structure
\begin{equation*}
\sigma(L_h) = \cup_{j \in \mathbb{N}} S_j,
\end{equation*}
where $S_j$ are the spectral bands associated with Bloch waves propagating inside the crystal. For each $\alpha \in Y^*$ the Bloch eigenvalues $\omega^2$ are of finite multiplicity and denoted by $\lambda_j(h, \alpha)$ with $\lambda_j(h, \alpha) \le \lambda_{j+1}(h, \alpha)$, $j \in \mathbb{N}$.

The band structure for the crystal is described by the family of dispersion relations
\begin{equation*}
    \omega^2 = \lambda_j(h, \alpha), \ j \in \mathbb{N}, \ \alpha \in Y^* 
\end{equation*}
and the spectral bands are given by the intervals
\begin{equation*}
S_j = [\min_{\alpha \in Y^*} \lambda_j(h, \alpha), \max_{\alpha \in Y^*} \lambda_j(h, \alpha)]. 
\end{equation*}
The upper and lower band edges of $S_j$ are denoted by
\begin{equation*}
b_j = \max_{\alpha \in Y^*} \lambda_j(h, \alpha) \quad \text{and} \quad a_j = \min_{\alpha \in Y^*} \lambda_j(h, \alpha) \end{equation*}
respectively. The band gaps are frequency intervals $\omega_- < \omega < \omega_+$ for which no waves propagate inside the crystal, i.e.,
\begin{equation*}
b_j < \omega_-^2 < \omega_+^2 < a_{j+1}. 
\end{equation*}

The corresponding eigenpair $(\lambda(h, \alpha), u_h)$ satisfies
\begin{equation*}
    L_h\,u_h=\lambda_j(h, \alpha)\,u_h.
\end{equation*}
We analyze the behavior of $\lambda_j(h, \alpha)$ for $|z|<U$, where the radius $U$ is explicitly determined by the inclusion geometry. Restricting back to the real axis $z=d^{2}\in \mathbb{R}$, further analysis establishes the existence of band gaps and pass bands at sufficiently high frequencies. These spectral features are computed explicitly and expressed in terms of the shape and placement of the physical inclusion within the $d$-periodic metamaterial crystal. 

\begin{figure}[h]
\centering
\includegraphics[scale=0.8]{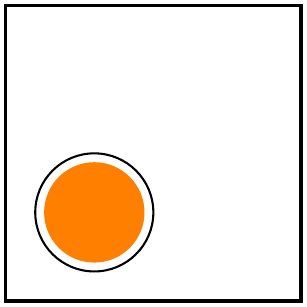}
\caption{Buffered disk inclusion.}
\label{inclusionsshaded}
\end{figure}

To illustrate the ideas, we now consider both Neumann and Bloch spectra for crystals for which each period cell contains a disk $D$ of radius $a$, see Figure~\ref{inclusionsshaded}. We suppose that the smallest distance that separates the disk from the boundary of the period cell is $t_d > 0$. The buffer layer $R$ is an annulus  with inner radius $a$ and outer radius $b = a+t$ where $t \le t_d/2$ chosen so that the buffered disk lies within the period cell. The benefit of considering this specific class of geometries is that the lower bound $\mu^* > -\frac{1}{2}$ on the quasistatic resonances can be explicitly stated, allowing for clear formulas for the radius of convergence \eqref{radiusofcva} in terms of $a$, $b$ and the zeros of the Bessel functions.

For this case, a suitable constant $\theta$ is computed in \cite{brunocomp91} and is given by
\begin{equation}\label{spectheta}
\theta = \frac{b^2-a^2}{b^2+a^2}. 
\end{equation}
Since $a < b$, we have $0 < \theta < 1$. 

We introduce the Dirichlet spectrum associated with the Laplace operator $-\Delta_D$ on the inclusion. We consider the part of the spectrum associated with eigenfunctions having nonzero average over the disk. These simple eigenvalues are denoted by
\begin{equation}
    \delta_{0j}^* = \left(\frac{\eta_{0j}}{a}\right)^2,
    \label{simplespecnonzero}
\end{equation}
where $\eta_{0j}$ are the zeros of the Bessel function of order zero $J_0$. The associated eigenfunctions are rotationally symmetric and given by the normalized Bessel functions of order zero:
\begin{equation*}
    u_{0j} = J_0(r\eta_{0j}/a) / (a\sqrt{\pi}J_1(\eta_{0j})).
\end{equation*}
The Dirichlet eigenvalues associated with mean zero eigenfunctions are of geometric multiplicity two and are denoted by $\nu_{nj} = (\eta_{nj}/a)^2$, where $\eta_{nj}$ is the $j$th zero of the $n$th Bessel function $J_n$, $1 \le n$. Next, we introduce the roots $\nu_{0j}$ of the spectral function
\begin{equation*}
    S(\nu) = \nu \sum_{j \in \mathbb{N}} \frac{a_{0j}^2}{\nu - \delta_{0j}^*} - 1, 
\end{equation*}
where $a_{0j} = \int_D u_{0j} \, dx$ are the nonzero averages of the rotationally symmetric normalized eigenfunctions $u_{0j}$.

We write 
\begin{equation}
    \sigma_N = \left\{ \bigcup_{j \in \mathbb{N}} \eta\nu_{0j} \right\} \cup \left\{ \bigcup_{(n,j) \in \mathbb{N}^2} \eta\nu_{nj} \right\}.\label{NeuSpecDisk} 
\end{equation}
and the Dirichlet spectrum $\sigma(-\Delta_D)$ given by
\begin{equation*}
    \sigma(-\Delta_D) = \left\{ \bigcup_{j \in \mathbb{N}} \delta_{0j}^* \right\} \cup \left\{ \bigcup_{(n,j) \in \mathbb{N}^2} \eta\nu_{nj} \right\}. 
\end{equation*}

We now provide an explicit condition on $d$ that is sufficient to open a band gap in the vicinity of $\eta\delta_{0j}^*$ together with explicit formulas describing its location and bandwidth. 

\begin{theorem} (Opening a band gap) Given $\delta_{0j}^*$ (see \eqref{simplespecnonzero}) define the set $\sigma_N^+$ to be elements $\eta\nu \in \sigma_N$ (see \eqref{NeuSpecDisk}) for which $\eta\nu >\eta\delta_{0j}^*$. The element in $\sigma_N^+$ closest to $\eta\delta_{0j}^*$ is denoted by $\eta\nu_{j+1}$. Set $s_j'$ according to
\begin{equation*}
    s_j'=\eta^{-1} s_j = \frac{1}{2} \min \left\{ |\eta^{-1}\nu_{j+1}^{-1} - \eta^{-1}\nu^{-1}|; \ \eta\nu \in \sigma_N \right\}. 
\end{equation*}
We define $\bar{r}_{j,\eta}$ as
\begin{equation}
    \bar{r}_{j,\eta} = \frac{\pi^2 \eta^{-1}s_j (b^2 - a^2)}{2(b^2 + a^2) + \pi^2 s_j (b^2 + 3a^2)}. \label{radiusdiskcase}
\end{equation}
Then, for $\tilde{L}_d$ given by \eqref{tildaLd}, one has the band gap
\begin{equation*}
    \sigma(\tilde{L}_d) \cap \left( \eta\delta_{0j}^*, \eta\nu_{j+1} \left( 1 - \frac{\eta\nu_{j+1} s_j}{\frac{\bar{r}_{j,1}}{d^2} - \eta} \right) \right) = \emptyset 
\end{equation*}
if
\begin{equation}
    d <  \sqrt{(\bar{r}_{j,1})|\eta|^{-1} \left( 1 + \frac{s_j \nu_{j+1}}{1 - \frac{\delta_{0j}^*}{\nu_{j+1}}} \right)^{-1}}. \label{boundkdiskcase}
\end{equation}
\end{theorem}

Next, we provide an explicit condition on $d$ sufficient for the persistence of a spectral band together with explicit formulas that describe its location and bandwidth. 

\begin{theorem}(Persistence of passbands). Given $\delta_{0j}^*$ (see \eqref{simplespecnonzero}) define the set $\sigma_N^-$ to be elements $\eta\nu \in \sigma_N$ (see \eqref{NeuSpecDisk}) for which $\eta\nu < \eta\delta_{0j}^*$. The element in $\sigma_N^-$ closest to $\eta\delta_{0j}^*$ is denoted by $\eta\nu_j$. Set $s_j'$ according to
\begin{equation*}
    s_j'=\eta^{-1} s_j = \frac{1}{2} \min \left\{ |\eta^{-1}(\delta_{j0}^*)^{-1} - \eta^{-1}\delta^{-1}|; \ \delta \in \sigma(-\Delta_D) \right\}.
\end{equation*}

Define $\underline{r}_{j,\eta}$ to be
\begin{equation*}
  \underline{r}_{j,\eta} = \frac{2\pi^2 \eta^{-1}s_j (b^2 - a^2)}{2(b^2 + a^2) + 2\pi^2 s_j (b^2 + 3a^2)}.
\end{equation*}

Then one has a passband in the vicinity of $d^2\eta\delta_{0j}^*$ and
\begin{equation*}
  \sigma(\tilde{L}_d) \supset \left[\eta\nu_j, \eta\delta_{0j}^* \left( 1 - \frac{\eta\delta_{0j}^* s_j}{\frac{\underline{r}_{j,1}}{d^2} - \eta} \right) \right]
\end{equation*}
if
\begin{equation*}
    d < \sqrt{(\underline{r}_{j,1})|\eta|^{-1} \left( 1 + \frac{s_j \delta_{0j}^*}{1 - \frac{\nu_j}{\delta_{0j}^*}} \right)^{-1}}.
\end{equation*}
\end{theorem}

\section{Diagonalization of $T^\alpha_h$, Power Series Representation for Operator and Eigenvalues, and Bounds for NP Eigenvalues}
\label{s:PowerSeriesRep}

The first step used in our method to find high frequency band gaps is to explicitly represent the operator $T^\alpha_h$ associated with the bilinear form $B_{\eta,h}(u,w)$ \eqref{Betah}. 

Let $L^2_\#(\alpha, Y)$  denote the space of all $\alpha$-quasi-periodic complex-valued functions in $L^2_{loc}(\mathbb{R}^d)$. For $\alpha \neq 0$, the eigenfunctions $u_d$ of \eqref{Helm-d3} lie in
$$H^1_\#(\alpha, Y) = \{ u \in H^1_{loc}(\mathbb{R}^d) : u \text{ is } \alpha\text{-quasi-periodic} \},$$
which becomes a Hilbert space under the inner product \begin{equation}
\label{H1innerprod}
\langle u, w \rangle = \int_Y \nabla u \cdot \nabla \overline{w}\,dy.
\end{equation} 
When $\alpha = 0$, the pair $u_d(y) = 1$, $\omega^2 = 0$ is a solution to \eqref{Helm-d3}. For this case, the remaining eigenfunctions associated with nonzero eigenvalues are orthogonal to constants in the $L^2(Y)$-inner product and belong to the set of square-integrable $Y$-periodic functions with zero average, denoted by $L^2_\#(0,Y)$. They also belong to the space
$$H^1_\#(0, Y) = \{ u \in H^1_{loc}(\mathbb{R}^d) : u \text{ is periodic}, \int_Y u \, dy = 0\}.$$
The space $H^1_\#(0, Y)$ is also a Hilbert space with inner product defined by \eqref{H1innerprod}.

For any $h \in \mathbb{C}$, we define the sesquilinear form $B_h: H^1_\#(\alpha, Y) \times H^1_\#(\alpha, Y) \rightarrow \mathbb{C}$ by
\begin{equation*}
%\label{Bhdefined}
    B_h(u,w) = h\int_{Y\setminus D} \nabla u(y) \cdot \nabla \overline{w(y)}\,dy + \eta\int_D \nabla u(y) \cdot \nabla \overline{w(y)}\,dy.
\end{equation*}

By the Riesz representation theorem, there exists a linear operator $T_h^\alpha : H^1_\#(\alpha, Y) \rightarrow H^1_\#(\alpha, Y)$ such that
\begin{equation}
\label{Thdefined}
    \langle T_h^\alpha u, w \rangle = B_h(u,w)
\end{equation}
for any $u,w \in H^1_\#(\alpha, Y)$.

As shown in \cite{RobertRobert1}, the space admits the orthogonal decomposition
\begin{align}
\label{OrthogonalTriple}
    H^1_\#(\alpha, Y) = W_{1,\alpha} \bigoplus W_{2,\alpha} \bigoplus W_{3,\alpha}
\end{align}
with subspaces characterized by vanishing gradients in $D$ or $Y\setminus D$, and $W_{3,\alpha}$ the orthogonal complement when $\alpha \neq 0$, and 
\begin{align}
     W_{1,0} & = \left\{ u \in H^1_\#(0, Y) \, : \, \nabla u = 0 \text{ in } D \right\} \\
    W_{2, 0} & = \left\{ u = \tilde{u} - \left( \frac{1}{|Y|} \int_D \tilde{u} dx \right) \mathbf{1}_Y \mid \tilde{u} \in \tilde{H}^1_0(D) \right\} \\
    W_{3, 0} & = \left ( W_{1,0} \bigoplus W_{2,0} \right )^\perp
\end{align}
when $\alpha = 0$.  
Here, $\tilde{H}^1_0(D)$ denotes the subspace of functions $H^1_0(D)$ extended by zero into $Y \setminus D$ and let $\mathbf{1}_Y$ be the indicator function of $Y$. 

In \cite{RobertRobert1}, the authors employed single- and double-layer potentials to obtain a resolution of the identity for the subspace $W_{3,\alpha}$ corresponding to a compact self-adjoint operator $T: W_{3,\alpha} \rightarrow W_{3,\alpha}$ which is given by the bilinear form
\begin{equation}
\label{Tdefined}
    \langle Tu,w \rangle = \frac{1}{2} \int_{Y\setminus D} \nabla u \cdot \nabla \overline{w} - \frac{1}{2} \int_D \nabla u \cdot \nabla \overline{w}.
\end{equation}
The spectrum $\sigma(T)$ of $T$ clearly lies between $-\frac{1}{2}$ and $\frac{1}{2}$, and coincides with that of the Neumann-Poincar\'e operator $(K_D^{-\alpha})^* : L^2(\partial D) \rightarrow L^2(\partial D)$ associated with the inclusion domain $D$, given by
\begin{equation*}
%\label{NeumannPoincaredefined}
    (K_D^{-\alpha})^*\phi[x] = \text{p. v. } \int_{\partial D} \frac{\partial G_\alpha(x,y)}{\partial n_x}\phi(y) d\sigma(y) \text{ , } x\in \partial D.
\end{equation*}
Here $G_\alpha(x,y)$ is the $\alpha$-quasi-periodic Green's function associated with the Laplacian $-\Delta$, that is, when $\alpha\neq0$, given by the formula
\begin{equation}
\label{Green'sFunctionquasi}
    G_\alpha(x,y) = -\sum \limits_{\mathbf{n} \in \mathbb{Z}^2} \frac{e^{i(2\pi \mathbf{n} + \alpha)\cdot(x-y)}}{|2\pi \mathbf{n} + \alpha|^2}
\end{equation}
and when $\alpha = 0$ is given by
\begin{equation}
\label{Green'sFunctionperiodic}
    G_0(x,y) = -\sum \limits_{\mathbf{n} \in \mathbb{Z}^2\setminus \{ \mathbf{0}\}} \frac{e^{i(2\pi \mathbf{n})\cdot(x-y)}}{|2\pi \mathbf{n}|^2}.
\end{equation}
Introduce the single layer potential
\begin{equation}
\label{Single layer}
    S_{D}[\rho(x)] = \int_{\partial D} G^\alpha(x,y)\rho(y)d\sigma(y),\,x\in Y
\end{equation}
with $S_{D}:H^{-1/2}(\partial D)\rightarrow W_{3\alpha}$ and the operator $S_D:$ is ibounded and nvertible \cite{RobertRobert1}, \cite{RobertRobert2}.
 
 The operator $T$ is compact on $W_{2,\alpha}$. To see this, we note by \cite{khavinsonputinarshapiro07} that $K_D^{-\alpha}$  is compact on $H^{-1/2}(\partial D)$ and  $T=S_{D}K_D^{-\alpha}S_{D}^{-1}$,  \cite{RobertRobert1}. So  its spectrum is discrete with one possible cluster point at zero, and we write $\sigma(T) = \{ \mu_i \}_{i \in \mathbb{N}}$, noting that $-\frac{1}{2} < \mu_i \leq \frac{1}{2}$.  Each $\mu_i$ has an associated spectral projection $P_{\mu_i}:W_{3,\alpha} \rightarrow W_{3,\alpha}$ which is finite dimensional (with the exception possibly of $\mu_i = 0$), yielding the following resolution of the identity for $W_{3,\alpha}$:
\begin{equation*}
%\label{W3Resolution}
I_{W_{3,\alpha}} = \sum \limits_{-\frac{1}{2} < \mu_i \leq \frac{1}{2}}
P_{\mu_i}
\end{equation*}
This resolution of the identity, together with the orthogonal triple decomposition of $H^1_\#(\alpha, Y)$ \eqref{OrthogonalTriple} and the simultaneous diagonalization of $T$ \eqref{Tdefined} and $T_h^\alpha$ \eqref{Thdefined} over $W_{3,\alpha}$, are used to establish the following theorem; see \cite{RobertRobert1}.

\begin{theorem}
    The operator $T_h^\alpha: H^1_\#(\alpha, Y) \rightarrow H^1_\#(\alpha, Y)$ associated with the sesquilinear form $B_h(u,w)$ (see \eqref{Thdefined}) is given by
\begin{equation}
\label{Tdiagonal}
\langle T_h^\alpha u,w \rangle = \langle h P_1u + \eta P_2u + \sum \limits_{-\frac{1}{2} < \mu_i \leq \frac{1}{2}} [h(1/2+ \mu_i) + \eta(1/2-\mu_i)]P_{\mu_i} u , w \rangle
\end{equation}
where $P_1$ and $P_2$ are the orthogonal projections onto $W_{1,\alpha}$ and $W_{2,\alpha}$.

Moreover, for fixed $\eta$, we see that $T_h^\alpha$ is bounded and invertible as long as
\begin{equation*}
%h(1/2+ \mu_i) + \eta(1/2-\mu_i)&\neq0\\
 z=d^2\neq\frac{(\mu_i+1/2)}{\eta(\mu_i-1/2)}.
\end{equation*}
In this case, \begin{equation*}
\langle (T_h^\alpha)^{-1} u,w \rangle = \langle h^{-1} P_1u + \eta^{-1} P_2u + \sum \limits_{-\frac{1}{2} < \mu_i \leq \frac{1}{2}} h^{-1}[(1/2+ \mu_i) + \eta h^{-1}(1/2-\mu_i)]^{-1}P_{\mu_i} u , w \rangle
\end{equation*}
or 
\begin{equation*}
\langle (T_z^\alpha)^{-1} u,w \rangle = \langle z P_1u + \eta^{-1} P_2u + \sum \limits_{-\frac{1}{2} < \mu_i \leq \frac{1}{2}} z[(1/2+ \mu_i) + \eta z(1/2-\mu_i)]^{-1}P_{\mu_i} u , w \rangle.
\end{equation*}
\end{theorem}

For future reference, we introduce the set $S$ of $z \in \mathbb{C}$ for which $T_z^\alpha$ is not invertible given by
\begin{equation}
\label{Sblocka}
S = \left\{\frac{(\mu_i + 1/2)}{\eta(\mu_i - 1/2)} \middle| -\frac{1}{2} < \mu_i < \frac{1}{2}\right\}.
\end{equation}

The Bloch spectral problem can then be written in operator form as
$$-\nabla \cdot (h(1-\chi_D) + \eta\chi_D)\nabla u = -\Delta_\alpha T_h^\alpha u = \nu(h)u,$$
with $u \in H_\#(\alpha, \mathbb{R}^d)$. We develop a power series in $z = \frac{1}{h}$ about $z=0$ for the spectrum of the family of operators given by $A^\alpha_\eta(z) = (T_z^\alpha)^{-1}(-\Delta_\alpha)^{-1}$, i.e.,
$$A^\alpha_\eta(z) = \left(z P_1 + \eta^{-1} P_2 + \sum \limits_{-\frac{1}{2} < \mu_i \leq \frac{1}{2}} z[(1/2+ \mu_i) + \eta z(1/2-\mu_i)]^{-1}P_{\mu_i}\right)(-\Delta_\alpha)^{-1}.$$
Here, the operator $(-\Delta_\alpha)^{-1}$ defined for all $\alpha \in Y^*$ is given by
$$(-\Delta_\alpha)^{-1}u(x) = -\int_Y G^\alpha(x,y)u(y)dy,$$
where $G^\alpha(x,y)$ is the Green's function for the quasi-periodic Laplace operator \eqref{Green'sFunctionquasi}.
The operator $A^\alpha_\eta(z)$ is self-adjoint for $h \in \mathbb{R}$ and is a family of compact bounded operators that take $L^2_\#(\alpha, Y)$ into itself.  Also, as in \cite{RobertRobert1}, we have the following lemma.
\begin{lemma} Spectrum of $A^\alpha_\eta(z)$
\begin{enumerate}
    \item $A^\alpha_\eta(z)$ is holomorphic on $\Omega_0 := \mathbb{C} \setminus S$, where $S$ is given by \eqref{Sblocka}, which lies on the negative real axis. The set $S$ consists of poles of $A^\alpha(z)$ with only one accumulation point at $z = -1/\eta$.
    \item For $z \in \Omega_0$, the spectrum of $A^\alpha_\eta(z)$, denoted by $\sigma(A^\alpha_\eta(z))$, consists of eigenvalues $\beta_j^\alpha(z)$ of finite multiplicity with a possible accumulation point at zero.
\end{enumerate}
\end{lemma}

From \cite{RobertRobert1}, we have the following theorem.

\begin{theorem} The Bloch eigenvalue problem \eqref{Helm-d3} for $-\nabla(h(1-\chi_D) + \eta\chi_D)\nabla$ defined on $H^1_\#(\alpha,Y)$ can be extended for values of the coupling constant $h$ with $z=h^{-1} \in \mathbb{C}\setminus S$, $S$ given by \eqref{Sblocka}, so that for each $\alpha \in Y^*$ the Bloch eigenvalues are of finite multiplicity and denoted by $\lambda_j(h,\alpha) = (\beta_j^\alpha(h^{-1}))^{-1}= (\beta_j^\alpha(z))^{-1}$, $j \in \mathbb{N}$ and the band structure
$$\lambda_j(h,\alpha) =\omega^2, \; j \in \mathbb{N}$$
extends to complex coupling constants $h$ with $z=h^{-1} \in \mathbb{C} \setminus S$.
\end{theorem}

\begin{figure}[ht]
\centering
\includegraphics{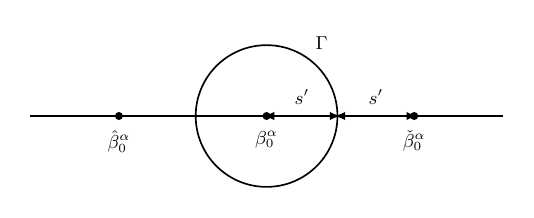}
\caption{Diagram of contour $\Gamma$}
\label{countourgraph1}
\end{figure}

Let $\beta_0^\alpha \in \sigma(A^\alpha_\eta(0))$ with spectral projection $P(0)$, and let $\Gamma$ be a closed contour in $\mathbb{C}$ enclosing $\beta_0^\alpha$ but no other element of $\sigma(A^\alpha_\eta(0))$, see Figure~\ref{countourgraph1}. 
The spectral projection associated with $\beta^\alpha(z) \in \sigma(A^\alpha_\eta(z))$ for $\beta^\alpha(z) \in \text{int}(\Gamma)$ is denoted by $P(z)$. We suppose that all elements $\beta_j^\alpha(z) \in \text{int}(\Gamma)$ converge to $\beta_0^\alpha$ as $z \to 0$. This collection is known as the eigenvalue group associated with $\beta_0^\alpha$. We write $M(z) = P(z)L^2_\#(\alpha,Y)$ and suppose for the moment that $\Gamma$ lies in the resolvent of $A^\alpha_\eta(z)$ and $\text{dim}(M(0)) = \text{dim}(M(z)) = m$. Define the resolvent of $A^\alpha_\eta(z)$ by
$$R(\zeta,z) = (A^\alpha_\eta(z) - \zeta)^{-1},$$
and, following the analysis in \cite{RobertRobert1}, we can obtain the following power series for the resolvent
\begin{equation}
\label{seriesresolvantblocka}
R(\zeta,z) = R(\zeta,0) + \sum_{n=1}^\infty z^n R_n(\zeta),
\end{equation}
which delivers a series representation formula for Bloch eigenvalues $\lambda_j(h,\alpha) = (\beta_j^\alpha(\frac{1}{h}))^{-1}$ for real $h$ in a neighborhood of $\infty$ given by
\begin{equation}
\label{seriesbetaiblocha}
\beta_j^\alpha(z) = \beta_0^\alpha + \sum_{n=1}^\infty z^n \beta_j^n(z).
\end{equation}

Also, we have $A_\eta^\alpha(0)= \eta^{-1} P_2\Delta^{-1}_{\alpha},$ which has the form
\begin{align*}
   \eta^{-1}\int_{D}\nabla u\cdot\nabla\overline{w}=\lambda\int_{Y}u\overline{w}. 
\end{align*}

The spectrum of $A_\eta^\alpha(0)$, for $\alpha\neq0$ is established in Theorem~\ref{th:blochspectanotzero-first} and for $\alpha=0$ is established in Theorem~\ref{th:blochspectazero-first}.

Define $s^\prime=\eta^{-1}s$ as the distance between $\Gamma$ and $\eta^{-1}\sigma(-\Delta_D^{-1})$, that is,
\begin{equation*}
    s^\prime =\eta^{-1}s = \inf_{\zeta \in \Gamma} \{\text{dist}(\zeta, \eta^{-1}\sigma(-\Delta_D^{-1}))\}.
\end{equation*}
 The power series \eqref{seriesresolvantblocka}, and thus \eqref{seriesbetaiblocha}, converges uniformly if 
 \begin{align}
 &|z| < r^*:=\frac{|\alpha|^2|z^*||\eta|^{-1}{s}}{\frac{1}{\frac{1}{2} - \mu^*} + |\alpha|^2{s}}\notag\\
 &\iff |d| < \sqrt{r^*}=\sqrt{\frac{|\alpha|^2|z^*||\eta|^{-1}{s}}{\frac{1}{\frac{1}{2} - \mu^*} + |\alpha|^2{s}}}\label{radiusofcva}.
\end{align}
As seen in the inequality $\eqref{radiusofcva}$, the radius of convergence for the power series expansions \eqref{seriesresolvantblocka} and \eqref{seriesbetaiblocha} depends explicitly on  $-1/2<\mu^*$, the lower bound for both the Neumann electrostatic resonance spectra and quasi-periodic resonance spectra, see \cite{RobertRobert1}, which hold for generic classes of inclusion domains $D$. The corresponding upper bound $z^*$ on $S$ \eqref{Sblocka} is written
\begin{equation}
\max_{i}\{z_i\} \leq \frac{\mu^* + 1/2}{\eta(\mu^* - 1/2)} = z^* < 0. \label{eq:4.4}
\end{equation}

We omit the proof of these statements as it follows the same derivation for the radii of convergence found in \cite{RobertRobert1, RobertRobert2}.

\section{Limit spectrum for periodic, quasi-periodic, and Neumann eigenvalues}
\label{s:limitspectra}

We begin by introducing the Neumann spectrum defined by the spectral problem on $Y$ given by
\begin{equation}
\label{Neumann}
-\nabla \cdot \left(\left( h(1-\chi_D(y))+\eta\chi_D(y)\right)\nabla u_h\right) = \nu(h) u_h(y), \quad y \in Y,
\end{equation}
together with the homogeneous Neumann boundary condition $\nabla u_h(y)\cdot \bf{n}=0$ on $\partial Y$, where $u_h$ is a function in $H^1(Y)$ with $\int_Y u_h\,dy = 0$. Here, $\chi_D$ denotes the indicator function of $D$ in the unit period. The Neumann eigenvalues for \eqref{Neumann} at fixed $ h > 1$ are written $\nu_j(h)$, $j=1,2,\dots$ and ordered according to min-max with $0 < \nu_1(h) \le \nu_2( h), \dots$. For future reference, we note for $h > 1$ and from the Min-Max principle that the Neumann spectrum provides a lower bound on the Bloch spectrum given by
\begin{equation}
\label{lbNtoD}
\nu_j( h) \le \lambda_j( h,\alpha), \text{ for } \alpha \in Y^*. 
\end{equation}

We begin by deriving a spectral representation formula for the inverse operator $(-\nabla \cdot ( h(1-\chi_D)+\eta\chi_D)\nabla)^{-1}$ associated with the Neumann spectrum. To proceed we complexify the problem and consider $ h \in \mathbb{C}$ noting that the divergence form operator $-\nabla \cdot ( h(1-\chi_D)+\eta\chi_D)\nabla$ is no longer uniformly elliptic. 

We denote the space of all square integrable complex-valued functions $u$ defined on $Y$ with $\int_Y u\,dy = 0$ by $L_0^2(Y)$ and the $L^2$-inner product over $Y$ is written
\[
(u,w) = \int_Y u\bar{w}\,dx. 
\]

The eigenfunctions $u_h$ for \eqref{Neumann} belong to the space
\[
\mathcal{H} = \{u \in H^1(Y) : \int_Y u\,dy = 0\}.
\]
The space $\mathcal{H}$ is a Hilbert space under the inner product
\[
\langle u,w \rangle = \int_Y \nabla u(y) \cdot \nabla \bar{w}(y)\,dy.
\]
The solutions to the Neumann spectral problem \eqref{Neumann} belong to the subspace $\mathcal{H}_N \subset \mathcal{H}$ defined by
\begin{equation*}
\mathcal{H}_N = \left\{ u \in \mathcal{H} : \frac{\partial u}{\partial \bf{n}} = 0 \text{ on } \partial Y \right\}. 
\end{equation*}
For any $h \in \mathbb{C}$, the weak formulation of the eigenvalue problem \eqref{Neumann} for $u_h$ and $ \nu(h)$ can be written as
\[
B_{h}(u_h,w) =  \nu(h) \langle u_h,w \rangle 
\]
for all $w \in \mathcal{H}_N$, where $B_{h} : \mathcal{H}_N \times \mathcal{H}_N \to \mathbb{C}$ is the sesquilinear form
\[
B_{h}(u,w) = { h} \int_{Y \setminus D} \nabla u(y) \cdot \nabla \bar{w}(y)\,dy + \eta\int_D \nabla u(y) \cdot \nabla \bar{w}(y)\,dy. 
\]
The linear operator $T_h: \mathcal{H}_N \to \mathcal{H}_N$ associated with $B_h$ is defined by
\[
\langle T_h u, w \rangle := B_h(u,w). 
\]
Proceeding as in \cite{RobertRobert1}, \cite{RobertRobert2}, one readily obtains a spectral decomposition of $T_h$ as in \eqref{Tdiagonal}.  From this decomposition and following the analysis again of \cite{RobertRobert1}, \cite{RobertRobert2}, one obtains a power series representation of the inverse operator $A(z)=T_h^{-1}(-\Delta_N)^{-1}$, where $-\Delta_N$ is the negative Laplacian acting on the space of mean-zero functions in $H^1(Y)$,
\begin{equation*}
    A(z) = \sum \limits_{n=0}^\infty z^nA_n \, \, ,
\end{equation*}
along with a power series expansion for the $i$th eigenvalue $\beta_i(z)$ of $A(z)$,
\begin{equation}\label{eigenexpneuman}
\beta_i(z) = \beta_0 + \sum_{n=1}^\infty z^n \beta_i^n(z).
\end{equation}
where 
\begin{equation}\label{betanu}
\beta_i(z) = \nu^{-1}(h)
\end{equation}
and $z=\frac{1}{h}$.  This power series will converge uniformly, provided that $z$ satisfies
\begin{equation}
\label{radiusneumanna}
    |z| < r^* = \frac{\lambda_N|z^*||\eta|^{-1} {s}}{\frac{1}{\frac{1}{2} - \mu^*} + \lambda_N{s}},
\end{equation} 
with $\lambda_N=\pi^2$, the first nonzero Neumann eigenvalue for the period $Y$. 

At this point, we identify the spectrum of the limiting operators $A(0)$ and $A^\alpha_\eta(0)$. Recall that the families of operators $A(z)$ and $A^\alpha_\eta(z)$ represent the inverse of the divergence form operator in the spectral problems \eqref{Neumann} and \eqref{Helm-d3}, indexed in $z = h^{-1}$, with Neumann and $\alpha$-quasi-periodic boundary conditions on $Y$, respectively. The limit spectra of these operators play a crucial role in establishing the upper and lower bounds of the band gaps and band edges identified in Theorem~\ref{openingbandgapbthm-first}.
Using the representation
$$A(z) =\left(z P_1 + \eta^{-1} P_2 + \sum \limits_{-\frac{1}{2} < \mu_i \leq \frac{1}{2}} z[(1/2+ \mu_i) + \eta z(1/2-\mu_i)]^{-1}P_{\mu_i}\right)(-\Delta_N)^{-1},$$
we see that
$$A(0) = \eta^{-1}P_2(-\Delta_N)^{-1},$$
and denote the spectrum of $A(0)$ by $\sigma(A(0))$.
To begin, consider the Dirichlet eigenvalues of the Laplace operator on $D$ associated with the spectral problem $-\Delta w = \delta w$, $\psi \in H^1_0(D)$ and denote the spectrum by $\sigma(-\Delta_D)$. The subset of Dirichlet eigenvalues associated with eigenfunctions that have zero mean over $D$ is denoted by $\{\delta_j^\prime\}_{j \in \mathbb{N}}$ and the set of Dirichlet eigenvalues associated with nonzero mean eigenfunctions is denoted by $\{\delta_j^*\}_{j \in \mathbb{N}}$. Next, we introduce the sequence of numbers $\{\nu_j\}_{j \in \mathbb{N}}$ given by the positive roots $\nu$ of the spectral function $S(\nu)$ defined by
\begin{equation}S(\nu) = \nu \sum_{i \in \mathbb{N}} \frac{a_i^2}{\nu - \delta_i^*} - 1,\label{spectfunct-first}
\end{equation}
where the coefficients $a_j = |\int_D \psi_j \, dx|$ are integrals of eigenfunctions $\psi_j$ corresponding to the Dirichlet eigenvalues $\delta_j$. The explicit characterization of $\sigma(A(0))$ is given by the following theorem; see \cite{RobertRobert1}.

\begin{theorem}
\label{spectrumNeumanb1}
   $\sigma(A(0)) = \{{\eta^{-1}\delta_j^\prime}^{-1}\}_{j \in \mathbb{N}} \cup \{\eta^{-1}\nu_j^{-1}\}_{j \in \mathbb{N}}$, where the families $\{{\delta_j^\prime}^{-1}\}_{j \in \mathbb{N}}$ and $\{\nu_j^{-1}\}_{j \in \mathbb{N}}$ are identified in the above paragraph. 
\end{theorem} 

We now recover the spectrum of the limiting operator $A^\alpha_\eta(0)$ when $\alpha \in Y^*$. Using the representation
$$A^\alpha_\eta(z) = (z P_1 + \eta^{-1} P_2 + \sum \limits_{-\frac{1}{2} < \mu_i \leq \frac{1}{2}} z[(1/2+ \mu_i) + \eta z(1/2-\mu_i)]^{-1}P_{\mu_i})(-\Delta_\alpha)^{-1},$$
we see that
$$A^\alpha(0) =\eta^{-1} P_2(-\Delta_\alpha)^{-1}.$$
The following theorem from \cite{RobertRobert1} provides the explicit characterization of the spectrum $\sigma(A^\alpha(0))$ for $\alpha \neq 0$.

\begin{theorem} \label{th:blochspectanotzero-first}Limit spectrum for quasi-periodic problem
$$\sigma(A^\alpha(0)) = \{{\eta^{-1}\delta_j^\prime}^{-1}\}_{j \in \mathbb{N}} \cup \{\eta^{-1}\delta_j^{*-1}\}_{j \in \mathbb{N}} = \eta^{-1}\sigma(-\Delta_D^{-1}),$$ where the families $\{{\delta_j^\prime}^{-1}\}_{j \in \mathbb{N}}$ and $\{\delta_j^{*-1}\}_{j \in \mathbb{N}}$ are identified in the paragraph above equation \eqref{spectfunct-first}.
\end{theorem}
To conclude, we recover the limit spectrum $\sigma(A^0(0))$ for the periodic problem, see \cite{RobertRobert1}.

\begin{theorem}\label{th:blochspectazero-first}$\sigma(A^0(0)) = \{\eta^{-1}\delta_j^{-1}\}_{j \in \mathbb{N}} \cup \{\eta^{-1}\nu_j^{-1}\}_{j \in \mathbb{N}}$, where the families $\{\delta_j^{-1}\}_{j \in \mathbb{N}}$ and $\{\nu_j^{-1}\}_{j \in \mathbb{N}}$ are identified in the paragraph above equation \eqref{spectfunct-first}.
\end{theorem}

\section{Upper bound on wavelength for convergence  of power series}
\label{s:wavelengthbound}

Given the relation $z=d^2 = \frac{| \alpha |^2\lambda^2}{4\pi^2}$ between the parameter $z$ and the quasi-momentum $\alpha$, and the inequality $\eqref{radiusofcva}$ that $|z|$ must satisfy in order for the power series $\eqref{seriesbetaiblocha}$ to converge uniformly, it follows that there is a bound on $|\alpha|$, and hence the wavelength $\lambda$, for which we can actually expand the Bloch eigenvalues in a uniformly convergent power series.  We determine that bound now.

Since $z=d^2 = \frac{| \alpha |^2\lambda^2}{4\pi^2}$, the inequality $|z| < r^*$ can be restated as follows.

\begin{align*}
    \frac{| \alpha |^2\lambda^2}{4\pi^2} &< \frac{|\alpha|^2|z^*||\eta|^{-1}{s}}{\frac{1}{\frac{1}{2} - \mu^*} + |\alpha|^2{s}}\\
    \iff 
    \lambda^2\left(1+ |\alpha|^2{s}\left(\frac{1}{2} - \mu^*\right) \right)  &< {4\pi^2}|\eta|^{-1}{s}\left(\frac{1}{2} + \mu^*\right) \\
    \iff 
    \lambda^2 |\alpha|^2{s}\left(\frac{1}{2} - \mu^*\right)   &< {4\pi^2}|\eta|^{-1}{s}\left(\frac{1}{2} + \mu^*\right)-\lambda^2  \\
   \iff|\alpha|^2 \, \, & < \, \, \frac{{4\pi^2}|\eta|^{-1}{s}\left(\frac{1}{2} + \mu^*\right)-\lambda^2}{\lambda^2 {s}\left(\frac{1}{2} - \mu^*\right)}. 
\end{align*}
 In order for $|\alpha|^2$ to have a  positive upper bound, we require
\begin{equation*}
\lambda^2 < {4\pi^2}|\eta|^{-1}{s}\left(\frac{1}{2} + \mu^*\right),
\end{equation*}
and since $k = 2\pi/\lambda$ delivers the equivalent result for the wave vector
\begin{equation*}
k > \sqrt{ \frac{|\eta|}{\left(\frac{1}{2} + \mu^*\right){s}}}.
\end{equation*}
This provides an upper bound on the wavelengths for which this analysis is valid, in terms of the Neumann-Poincar\'e spectrum of the inclusion geometry. 

\section{Band gaps}
\label{s:Bandgaps}
In this section, we combine the power series convergence results and the corresponding error analysis to provide a rigorous proof of the existence of band gaps in the spectrum of the high-contrast metamaterial. The proof relies on showing that the gaps present in the limit spectrum (as $z \to 0$) persist for sufficiently small but finite $z$.

\subsection{Uniform Error Estimates}

We begin by stating the uniform error bound required to control the dispersion of the Bloch bands across the Brillouin zone.

\begin{lemma}
Let $z$ be the high-contrast parameter and $R$ be the radius of convergence determined by the Neumann-Poincaré spectrum and the geometry of the inclusion $D$. For $z$ such that $|z| < \gamma R$ with $0 < \gamma < 1$, there exists a constant $C$, independent of the quasiperiodicity $\alpha$, such that the $n$-th Bloch eigenvalue $\lambda_n(z, \alpha)$ satisfies:
\begin{equation*}
\left| \lambda_n(z, \alpha) - \sum_{m=0}^{N} z^m \lambda_n^{m}(\alpha) \right| \leq C \frac{|z|^{N+1}}{1 - \frac{|z|}{R}}.
\end{equation*}\label{katoestimate}
\end{lemma}

\subsection{Proof of Gap Persistence}\label{sec:gapPersistence}

To complete the band-gap analysis, we utilize the Neumann problem on the unit cell $Y$ developed in the preceding sections. This problem provides an $\alpha$-independent lower bound on the Bloch spectrum. 

We apply the characterization of $\sigma(A(0))$ given by Theorem~\ref{spectrumNeumanb1} together with \eqref{betanu} to recover the high-contrast limit of the Neumann spectrum given by
\begin{equation}
   \sigma_N = \{\eta\delta'_j\}_{j\in\mathbb{N}} \cup \{\eta\nu_j\}_{j\in\mathbb{N}}, \label{Neumanspectbcharact1}
\end{equation}
where $\delta'_j$ is the part of the Dirichlet spectra for $D$ associated with mean-zero eigenfunctions and $\nu_i$ are the roots of the spectral function \eqref{spectfunct-first}. In what follows, we do not distinguish between the two component parts of the spectrum and write elements of $\sigma_N$ as $\eta\nu_j$, $j \in \mathbb{N}$. The Dirichlet spectrum is given by $\sigma(-\Delta_D) = \{\delta^\prime_j\}_{j\in\mathbb{N}} \cup \{\delta^*_j\}_{j\in\mathbb{N}}$ where $\delta^*_j$ are Dirichlet eigenvalues associated with nonzero mean eigenfunctions. The relation between $\sigma_N$ and $\eta\sigma(-\Delta_D)$ is given by the following theorem; see \cite{hempel2000spectral}.

\begin{theorem}\label{thm:interlacing}
Strict interlacing of the spectra.  If $\delta^*_j \in \sigma(-\Delta_D)$ is simple, then there exist adjacent elements $\eta\nu_j < \eta\nu_{j+1}$ ordered by min-max belonging to $\sigma_N$ such that
\begin{equation}
    \eta\nu_j < \eta\delta^*_j < \eta\nu_{j+1}.  \label{interlacing-first}
\end{equation}
\end{theorem}

Theorem~\ref{thm:interlacing} ensures the existence of a band gap for sufficiently large contrast $h$. In Theorem~\ref{openingbandgapbthm-first}, we provide an explicit condition on the contrast $h$ that guarantees the opening of a band gap in a neighborhood of $\eta \delta^*_j$, along with explicit formulas describing its location and bandwidth. Before stating the theorem, we introduce the following definition:
\begin{definition}\label{def:ptheta} The class of crystal geometries characterized by inclusions $D$ such that, for all $u \in W_{3,\alpha}$,
\begin{equation}
\|\nabla u\|_{L^2(Y \setminus D)}^2 \ge \theta \|\nabla u\|_{L^2(D)}^2,
\label{eq:8.2}
\end{equation} holds for a fixed positive value of $\theta$, is denoted by $P_\theta$.
\end{definition}

\begin{theorem}\label{openingbandgapbthm-first}
{\bf Opening a band gap. } Consider any crystal geometry belonging to the class $P_\theta$ given in Definition~\ref{def:ptheta}. Suppose $\eta\delta^*_j \in \eta\sigma(-\Delta_D)$ is simple, so that $\eta\delta^*_j < \eta\nu_{j+1}$ as in Theorem~\ref{thm:interlacing}. Set
$$ s_j'=\frac{1}{2} \mathrm{dist}(\{\eta^{-1}\nu^{-1}_{j+1}\}, \sigma_N \setminus \{\eta^{-1}\nu^{-1}_{j+1}\}) $$
with $s_j'=\eta^{-1}s_j$, and where $\sigma_N$ is given by \eqref{Neumanspectbcharact1}, and
\begin{equation}
\label{radiusneumannbbar-first}
    \bar{r}_{j,\eta}  = \frac{\pi^2|z^*||\eta|^{-1} {s_j}}{\frac{1}{\frac{1}{2} - \mu^*} + \pi^2{s_j}},
\end{equation} 
Then one has the band gap 
\begin{equation}\label{mr:bandgap}
    \sigma(L_h) \cap \left(\eta\delta^*_j, \eta\nu_{j+1}\left(1 - \frac{\eta\nu_{j+1} s_j}{h\bar{r}_{j,1} - \eta}\right)\right) = \emptyset
\end{equation}
if
\begin{equation}\label{mr:valuesh}
  h > \bar{h}_j = (\bar{r}_{j,1})^{-1}|\eta| \left(1 + \frac{s_j \nu_{j+1}}{1 - \frac{\delta^*_j}{\nu_{j+1}}}\right). 
\end{equation}
\end{theorem}
\begin{proof} Consider $\eta\nu_{j+1}$ of multiplicity $m < \infty$. By results on the separation of Neumann spectra and radius of convergence \eqref{radiusneumannbbar-first} for \eqref{eigenexpneuman}, we have that any element $\beta_i(z)$ in the eigenvalue group has series representation given by
\begin{equation*}
    \beta_i(z) = \eta^{-1}\nu_{j+1}^{-1} + \sum_{n=1}^\infty z^n \beta_i^{n}(z). 
\end{equation*} 
For $z$ inside the interval $-\bar{r}_{j,\eta} < z < \bar{r}_{j,\eta}$, for $h > (\bar{r}_{j,\eta})^{-1}$, let $\nu_{j+1}(h)$ be the minimum of $1/\beta_i(h)$, with $\beta_i(h)$ in the eigenvalue group, and we have $\nu_{j+1}(h) \to \eta\nu_{j+1}$ for $h \to \infty$. Since $\nu_{j+1}(h)$ is increasing with $h$, we conclude
\begin{equation*}
    \nu_{j+1}(h) \le \eta\nu_{j+1}. 
\end{equation*} 
We take a min-max order for $\eta\sigma(-\Delta_D)$ and suppose that the $j^{\text{th}}$-Dirichlet eigenvalue corresponds to an eigenfunction of nonzero mean. The eigenvalue is denoted by $\eta\delta^*_j$. The min-max principle together with the monotonicity of eigenvalues with respect to increasing $h$ delivers the inequality between the Bloch eigenvalues and $\eta\sigma(-\Delta_D)$:
\begin{equation}
\label{lambdajhaleqdjstar-first}
    \lambda_j(h, \alpha) \le \eta\delta^*_j, \text{ for } \alpha \in Y^*, \text{ and } h > 0. 
\end{equation} 
\begin{figure}[ht]
\centering
\includegraphics{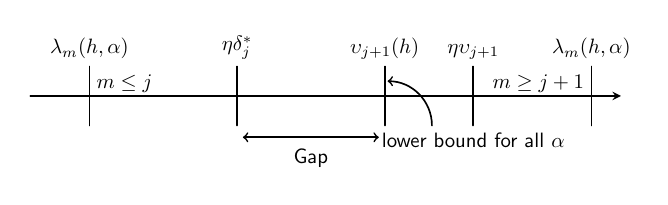}
\caption{Diagram for proof of Theorem~\ref{openingbandgapbthm-first}}
\label{bandgapopening-first}
\end{figure}

Application of \eqref{lbNtoD} and \eqref{lambdajhaleqdjstar-first} gives
\begin{equation*}
    \lambda_m(h, \alpha) \le \eta\delta^*_j, \text{ for } m \le j \quad \text{and} \quad \nu_{j+1}(h) \le \lambda_m(h, \alpha), \text{ for } j+1 \le m,  
\end{equation*} 
for every $\alpha \in Y^*$ and it is clear that a band gap opens in the Bloch spectrum when $\eta\delta^*_j < \nu_{j+1}(h)$ or equivalently when
\begin{equation}
\label{triangineqbandgapa}
    |\nu_{j+1}(h) - \eta\nu_{j+1}| < |\eta\delta^*_j - \eta\nu_{j+1}|=|\eta||\delta^*_j - \nu_{j+1}|, 
\end{equation} 
see Figure~\ref{bandgapopening-first}.

We apply Lemma~\ref{katoestimate} to get
\begin{align*}
|\nu_{j+1}(h) - \eta\nu_{j+1}| &< \frac{|\eta|^2\nu_{j+1}^2|\eta|^{-1} s_j}{h\bar{r}_{j,\eta} - 1}\\
&= \frac{|\eta|\nu_{j+1}^2 s_j}{h|\eta|^{-1}\bar{r}_{j,1} - 1}
\end{align*}
and \eqref{mr:bandgap} follows for all $h$ that satisfy
$$ \frac{\nu_{j+1}^2 s_j}{h|\eta|^{-1}\bar{r}_{j,1} - 1} < |\delta^*_j - \nu_{j+1}|, $$
which are given by \eqref{mr:valuesh}.
\end{proof}

This confirms the existence of gaps in the high-frequency spectrum, with the size and location determined by the contrast and the geometry of the inclusion phase $D$.

\begin{theorem}Persistence of passbands.

Consider any crystal geometry that belongs to the class $P_\theta$ defined by Definition~\ref{def:ptheta}. Suppose $\eta\delta_j^*$ is simple, so that $\eta\nu_j < \eta\delta_j^*$ as in Theorem~\ref{thm:interlacing}. Set
\[
s_j' = \frac{1}{2} \text{dist} \left( \{\eta^{-1}(\delta_j^*)^{-1}\}, \eta^{-1}\sigma(-\Delta_D) \setminus \{\eta^{-1}(\delta_j^*)^{-1}\} \right),
\]
where $s_j'=\eta^{-1}s_j$
and
\begin{equation}
\underline{r}_{j,\eta} = \frac{d \pi^2 |\eta|^{-1}s_j |z^*|}{\frac{1}{1/2 - \mu^*} + d \pi^2 s_j}, \text{ for } d = 2, 3.\label{underr}
\end{equation}
Then one has a passband in the vicinity of $\eta\delta_j^*$ and
\begin{equation}
\sigma(L_h) \supset \left[ \eta\nu_j, \eta\delta_j^* \left( 1 - \frac{\eta\delta_j^* s_j}{h \underline{r}_{j,1} - \eta} \right) \right] \label{sigmalhper}
\end{equation}
if
\begin{equation}
h > \underline{h}_j = (\underline{r}_{j,1})^{-1} |\eta|\left( 1 + \frac{s_j \delta_j^*}{1 - \frac{\nu_j}{\delta_j^*}} \right). \label{sigmalhperh}
\end{equation}
\end{theorem}
\begin{proof}
From the min-max formulation we see that $\lambda_j(h,0)$ and $\lambda_j(h,\alpha)$ are increasing with $h$ and from \cite{RobertRobert1} we have
\begin{equation*}
\lim_{h \to \infty} \lambda_j(k,0) = \eta\nu_j \text{ hence } \lambda_j(h,0) \leq \eta\nu_j, \text{ for } h > 1 
\end{equation*}
and
\begin{equation*}
\lim_{h \to \infty} \lambda_j(h,\alpha) = \eta\delta_j^* \text{ hence } \lambda_j(h,\alpha) \leq \eta\delta_j^*, \text{ for } h > 1. 
\end{equation*}
With this in mind, observe that if $|\eta\delta_j^* - \lambda_j(h,\alpha)| < |\eta\delta_j^* - \eta\nu_j|$ then
\begin{equation}
\sigma(L_h) \supset [\eta\nu_j, \lambda_j(h,\alpha)].\label{persbandstr} 
\end{equation}
To proceed, we estimate the difference $|\lambda_j(h,\alpha) - \eta\delta_j^*|$. Set 
\[
s_j' = \frac{1}{2} \text{dist} \left( \{\eta^{-1}(\delta_j^*)^{-1}\}, \eta^{-1}\sigma(-\Delta_D) \setminus \{\eta^{-1}(\delta_j^*)^{-1}\} \right),
\]
where $s_j'=\eta^{-1}s_j$ and let $r_j^*=r_{j,\eta}^*$ be given by \eqref{radiusofcva}. Noting that $\eta\delta_j^*$ is simple, $\beta_j^\alpha(z)$ has series representation given by
\begin{equation*}
\beta_j^\alpha(z) = \eta^{-1}(\delta_j^*)^{-1} + \sum_{n=1}^{\infty} z^n \beta_j^n(z), 
\end{equation*}
for $z$ inside the interval $-r_j^* < z < r_j^*$. For $h > {r_j^*}^{-1}$ we have $\lambda_j(h,\alpha) = \frac{1}{\beta_j^\alpha(h^{-1})}$ and we apply Lemma~\ref{katoestimate} to get
\begin{align*}
|\eta\delta_j^* - \lambda_j(h,\alpha)|& < \frac{|\eta|^2(\delta_j^*)^2|\eta|^{-1} s_j}{hr_{j,\eta}^* - 1}\\
& = \frac{|\eta|(\delta_j^*)^2 s_j}{h|\eta|^{-1}r_{j,1}^* - 1}
\end{align*}
The persistence of band structure described by \eqref{persbandstr} follows for a fixed $\alpha \in Y^*$ for all $h$ that satisfy
\begin{align*}
&\frac{|\eta|(\delta_j^*)^2 s_j}{h|\eta|^{-1}r_{j,1}^* - 1} < |\eta\delta_j^* - \eta\nu_j| \\
&\implies \frac{(\delta_j^*)^2 s_j}{h|\eta|^{-1}r_{j,1}^* - 1} < |\delta_j^* - \nu_j|.
\end{align*}
We maximize $r_j^*$ over $\alpha \in [-\pi, \pi]^d$ to find that it is attained for $|\alpha|^2 = d \pi^2$ and the maximum is $r_j^* = \underline{r}_j$. For this choice, we recover the persistence of band structure described by \eqref{sigmalhper} and \eqref{sigmalhperh}. 
\end{proof}
\section{Conclusion}
\label{sec:conclusion}

In this paper, we have rigorously investigated the high-frequency Bloch spectrum for transverse electric (TE) modes in a two-dimensional, high-contrast periodic dielectric metamaterial. By departing from standard sub-wavelength homogenization paradigms and moving into a regime where the incident wavelength is commensurate with the microstructural lattice period $d$, we successfully captured multiple scattering and cell resonance features across the entire Brillouin zone. Utilizing an operator factorization technique coupled with a three-space orthogonal Sobolev decomposition, we showed that the Bloch eigenvalues $\omega^2$ can be analytically continued from $d=0$ and explicitly represented via a convergent power series in $z = d^2$. 

A core contribution of this work is the development of a geometric, explicit threshold $\overline{d}$ below which this power series expansion is guaranteed to remain valid. This threshold is explicitly mapped to physical microstructural parameters governed by the interaction between the Dirichlet spectra of the inclusion and the electrostatic resonances of the periodic matrix, the latter tracked via the spectrum of the associated Neumann-Poincar\'{e} operator. Furthermore, we provided a definitive proof for the opening of absolute band gaps for all scales $d < \overline{d}$, and applied this methodology to the case of buffered disk inclusions. 

%Note that both thresholds $\overline{h}_j$ and $\underline{h}_j$ depend explicitly on the crystal geometry through $a$ and $b$ and the Dirichlet spectrum of the disk. Density results on the distribution of zeros of Bessel functions show that the distance between adjacent eigenvalues $s_j$ approaches zero with increasing $j$. This, together with \eqref{radiusdiskcase} and \eqref{boundkdiskcase}, implies that the contrast sufficient to open gaps grows without bound as $j \to \infty$. 

\bibliographystyle{habbrv}
\bibliography{homogenisation}
%\bibliography{../Bib/homogenisation}
%\bibliographystyle{abbrv}

\end{document}